\documentclass[11pt]{amsart}
\usepackage{amsmath,amsthm,amsfonts,amssymb,mathrsfs}
\date{\today}

\usepackage{hyperref}

\newtheorem{theorem}{Theorem}

\newtheorem{proposition}[theorem]{Proposition}
\newtheorem{corollary}[theorem]{Corollary}
\theoremstyle{definition}
\newtheorem{example}[theorem]{Example}

\newtheorem{remark}[theorem]{Remark}
\newtheorem{lemma}[theorem]{Lemma}

\begin{document}

\title[On linearly ordered $H$-closed topological semilattices]
{On linearly ordered $H$-closed topological semilattices}

\author{Oleg~Gutik}
\address{Department of Mathematics, Ivan Franko Lviv National University,
Universytetska 1, Lviv, 79000, Ukraine}
\email{o\_\,gutik@franko.lviv.ua}
\author{Du\v{s}an~Repov\v{s}}
\address{Institute of Mathematics, Physics and Mechanics, and Faculty of Education,
University of Ljubljana, P.O.B. 2964, Ljubljana, 1001, Slovenia}
\email{dusan.repovs@guest.arnes.si}

\keywords{Topological semilattice, linearly ordered topological
semilattice, $H$-closed topological semilattice, absolutely
$H$-closed topological semilattice.}

\subjclass[2000]{06A12, 06F30, 22A15, 22A26, 54H12}

\begin{abstract}
We give a criterium when a linearly ordered
topological semilattice is $H$-closed. We also
prove that any linearly ordered
$H$-closed topological semilattice is absolutely $H$-closed and we
show that every linearly ordered semilattice is a dense
subsemilattice of an $H$-closed topological semilattice.
\end{abstract}

\maketitle

In this paper all topological spaces will be assumed to be
Hausdorff. We shall follow the terminology of~\cite{Birkhoff}--\cite{GHKLMS}. If $A$ is a subset of a topological space
$X$, then we shall
denote
  the closure of the
set $A$ in $X$
by $\operatorname{cl}_X(A)$. We shall denote the first infinite cardinal
by $\omega$.

A {\em semilattice} is a semigroup with a commutative idempotent
semigroup operation. If $S$ is a topological space equipped with a
continuous semigroup operation, then $S$ is called a {\em
topological semigroup}. A {\em topological semilattice} is a
topological semigroup which is algebraically a semilattice.

If $E$ is a semilattice, then the semilattice operation on $E$
determines the partial order $\leqslant$ on $E$: $$e\leqslant
f\quad\text{if and only if}\quad ef=fe=e.$$ This order is called
{\em natural}. An element $e$ of a semilattice $E$ is called {\em
minimal} ({\em maximal}) if $f\leqslant e$ ($e\leqslant f$)
implies $f=e$ for $f\in E$. For elements $e$ and $f$ of a
semilattice $E$ we write $e<f$ if $e\leqslant f$ and $e\neq f$. A
semilattice $E$ is said to be {\em linearly ordered} or a
\emph{chain} if the natural order on $E$ is linear.

Let $S$ be a semilattice and $e\in S$. We denote ${\downarrow}
e=\{ f\in S\mid f\leqslant e\}$ and ${\uparrow} e=\{ f\in S\mid
e\leqslant f\}$. If $S$ is a topological semilattice then
Propositions~VI.1.6(ii) and VI.1.14 of~\cite{GHKLMS} imply that
${\uparrow} e$ and ${\downarrow} e$ are closed subsets in $S$ for
any $e\in S$.

Let $E$ be a linearly ordered topological semilattice. Since
${\downarrow} e$ and ${\uparrow} e$ are closed for each $e\in E$,
it follows that the topology of $E$ refines the order topology.
Thus the principal objects whick
we shall
consider can
be alternatively viewed
as linearly ordered sets equipped with the semilattice operation
of taking minimum and equipped with a topology refining the order
topology for which the semilattice operation is continuous.

Let ${\mathscr S}$ be some class of topological semigroups. A
semigroup $S\in{\mathscr S}$ is called {\it $H$-closed in}
${\mathscr S}$, if $S$ is a closed subsemigroup of any topological
semigroup $T\in{\mathscr S}$ which contains $S$ both as a
subsemigroup and as a topological space. If ${\mathscr S}$
coincides with the class of all topological semigroups, then the
semigroup $S$ is called {\it $H$-closed}. The $H$-closed
topological semigroups were introduced by Stepp
\cite{St1}, 
they were called {\it maximal semigroups}. A
topological semigroup $S\in{\mathscr S}$ is called {\it absolutely
$H$-closed in the class} ${\mathscr S}$, if any continuous
homomorphic image of $S$ into $T\in{\mathscr S}$ is $H$-closed in
${\mathscr S}$. If ${\mathscr S}$ coincides with the class of all
topological semigroups, then the semigroup $S$ is called {\it
absolutely $H$-closed}.

An algebraic semigroup $S$ is
called {\it algebraically $h$-closed in} ${\mathscr S}$, if $S$ equipped 
with discrete topology ${\mathfrak{d}}$ is absolutely $H$-closed
in ${\mathscr S}$ and $(S,{\mathfrak{d}})\in{\mathscr S}$. If
${\mathscr S}$ coincides with the class of all topological
semigroups, then the semigroup $S$ is called {\it algebraically
$h$-closed}. Absolutely $H$-closed topological semigroups and
algebraically $h$-closed semigroups were introduced by Stepp
in~\cite{St2}, they were called {\it absolutely maximal}
and {\it algebraic maximal}, respectively. 
Gutik and Pavlyk \cite{GP2}
observed
that a topological semilattice
is [absolutely] $H$-closed if and only if it is [absolutely]
$H$-closed in the class of topological semilattices.

Stepp
\cite{St2} proved that a semilattice $E$ is
algebraically $h$-closed if and only if any maximal chain in $E$
is finite and  he asked the following
question: {\it Is every
$H$-closed
topological semilattice absolutely $H$-closed?} 
In the present
paper we
give a criterium when  a linearly ordered topological
semilattice is
$H$-closed. 
We
also
prove that every
linearly ordered $H$-closed
topological semilattice is absolutely $H$-closed, we
show that every
linearly ordered semilattice is a dense subsemilattice of an
$H$-closed topological semilattice, and we
give an example of a
linearly ordered $H$-closed locally compact topological
semilattice which does
not embed into a compact topological
semilattice.

Let $C$ be a maximal chain of a topological semilattice $E$. Then
$C=\displaystyle\bigcap_{e\in C}({\downarrow} e\cup{\uparrow} e)$,
and hence $C$ is a closed subsemilattice of $E$. Therefore we obtain the following:

\begin{lemma}\label{lemma1}
Let $L$ be a linearly ordered subsemilattice of a topological
semilattice $E$. Then $\operatorname{cl}_E(L)$ is a linearly
ordered subsemilattice of $E$.
\end{lemma}

A linearly ordered topological semilattice $E$ is called
\emph{complete} if every non-empty subset of
$S$ has $\inf$ and
$\sup$.

\begin{theorem}\label{theorem1}
A linearly ordered topological semilattice $E$ is $H$-closed if
and only if the following conditions hold:
\begin{itemize}
    \item[$(i)$] $E$ is complete;
    \item[$(ii)$] $x=\sup A$ for $A={\downarrow}A\setminus\{ x\}$
    implies $x\in\operatorname{cl}_EA$, whenever
    $A\neq\varnothing$; and
    \item[$(iii)$] $x=\inf B$ for $B={\uparrow}B\setminus\{ x\}$
    implies $x\in\operatorname{cl}_EB$, whenever
    $B\neq\varnothing$
\end{itemize}
\end{theorem}

\begin{proof}
$(\Leftarrow)$ Suppose to the contrary
that
there exists a linearly
ordered non-$H$-closed topological semilattice $E$ which satisfies
the conditions $(i)$, $(ii)$, and $(iii)$. Then by
Lemma~\ref{lemma1} there exists a linearly ordered topological
semilattice $T$ such that $E$ is a dense proper subsemilattice of
$T$. Let $x\in T\setminus E$. Condition $(i)$ implies that
$x\neq\sup E$ and $x\neq\inf E$, otherwise $E$ is not complete.
Therefore we have $\inf E<x<\sup E$. Let
$A(x)=(T\setminus{\uparrow}x)\cap E$ and
$B(x)=(T\setminus{\downarrow}x)\cap E$. Since $E$ is complete, we
have $\sup A(x)\in E$ and $\inf B(x)\in E$. Let $s=\sup A(x)$ and
$i=\inf B(x)$. We observe that $s<x<i$, otherwise, if $x<s$, then
$A={\downarrow}x\cap E=({\downarrow}s\setminus\{ s\})_E$ is a
closed subset in $E$, which contradicts the
condition $(ii)$, and if
$i<x$, then $B={\uparrow}x\cap E=({\uparrow}i\setminus\{ i\})_E$
is a closed subset of $E$, which
contradicts the
condition $(iii)$.
Then ${\uparrow}i$ and ${\downarrow}s$ are closed subsets of $T$
and $T={\downarrow}s\cup{\uparrow}i\cup\{ x\}$. Therefore $x$ is
an isolated point of $T$. This contradicts
the assumption
that
$x\in T\setminus E$ and $E$ is a dense subspace of $T$. This
contradiction implies that $E$ is an $H$-closed
topological semilattice.

$(\Rightarrow)$ At first we shall show that $\sup A\in E$ for
every infinite subset $S$ of $E$. Suppose to the contrary, that
there exists an infinite subset $A$ in $E$ such that $A$ has no
$\sup$ in $E$. Since $S$ is a linearly ordered semilattice, we
have that the set ${\downarrow}A$ also has no $\sup$ in $E$. We
consider two cases:
\begin{itemize}
    \item[$(a)$] $E\setminus{\downarrow}A\neq\varnothing$; and
    \item[$(b)$] $E\setminus{\downarrow}A=\varnothing$.
\end{itemize}

In case $(a)$  the set
$B=E\setminus{\downarrow}A$ has no $\inf$,
since $E$ is a linearly ordered semilattice, 
otherwise $\inf B=\sup
A\in E$.

Let $x\notin E$. We put $E^*=E\cup\{ x\}$. We extend the
semilattice operation from $E$ onto $E^*$ as follows:
\begin{equation*}
    x\cdot y=y\cdot x=
    \left\{%
\begin{array}{ll}
    x, & \hbox{if}~~y\in E\setminus{\downarrow}A; \\
    x, & \hbox{if}~~y=x; \\
    y, & \hbox{if}~~y\in{\downarrow}A .\\
\end{array}%
\right.
\end{equation*}
Obviously, the semilattice operation on $E^*$ determines a linear
order on $E^*$.

We define a topology $\tau^*$ on $E^*$ as follows. Let $\tau$ be
the topology on $E$. At any point $a\in E=E^*\setminus\{ x\}$
bases of topologies $\tau^*$ and $\tau$ coincide. We put
\begin{equation*}
 \mathscr{B}^\ast(x)=\left\{
 V_b^a(x)=E\setminus\left({\downarrow}b\cup{\uparrow}a\right)\mid
 a\in E\setminus{\downarrow}A, b\in {\downarrow}A\right\}.
\end{equation*}
Since the set $E\setminus{\downarrow}A$ has no $\inf$ and the set
${\downarrow}A$ has no $\sup$, the conditions (BP1)---(BP3) of
\cite{En} hold for the family $\mathscr{B}^*(x)$ and
$\mathscr{B}^*(x)$ is a base of a Hausdorff topology $\tau^*$ at
the point $x\in E^*$.

Let $c\in E\setminus{\downarrow}A$ and $d\in{\downarrow}A$. Then
there exist $a\in E\setminus{\downarrow}A$ and $b\in{\downarrow}A$
such that $d<b<x<a<c$. Then for any open neighbourhoods $V(c)$ and
$V(d)$ of the points $c$ and $d$, respectively, such that
$V(c)\subseteq E\setminus{\downarrow}a=E^*{\downarrow}a$ and
$V(d)\subseteq E\setminus{\uparrow}b=E^*\setminus{\uparrow}b$ we
have
\begin{equation*}
   V_b^a(x)\cdot V(d)\subseteq V(d) \qquad \mbox{ and } \qquad
   V_b^a(x)\cdot V(c)\subseteq V_b^a(x).
\end{equation*}
We
also
have $V_b^a(x)\cdot V_b^a(x)\subseteq V_b^a(x)$ for all
$a\in E\setminus{\downarrow}A$ and $b\in{\downarrow}A$. Therefore
$(E^*,\tau^*)$ is a Hausdorff topological semilattice which
contains $E$ as a dense non-closed subsemilattice. This
contradicts the assumption that $E$ is an $H$-closed topological
semilattice.

Consider case $(b)$. Let $y\notin E$. We put $E^\star=E\cup\{ y\}$
and extend the semilattice operation from $E$ onto $E^\star$ as
follows:
\begin{equation*}
    y\cdot s=s\cdot y=
    \left\{%
\begin{array}{ll}
    y, & \hbox{if~} s=y; \\
    s, & \hbox{if~} s\neq y.\\
\end{array}%
\right.
\end{equation*}
Obviously, the semilattice operation on $E^\star$ determines a
linear order on $E^\star$.

We define a topology $\tau^\star$ on $E^\star$ as follows. Let
$\tau$ be the topology on $E$. At any point $a\in
E=E^\star\setminus\{ x\}$ bases of topologies $\tau^\star$ and
$\tau$ coincide. We put
\begin{equation*}
 \mathscr{B}^\star(x)=\left\{ V_b(x)=E\setminus{\downarrow}b\mid
 b\in{\downarrow}A=E\right\}.
\end{equation*}
Since the set $E={\downarrow}E$ has no $\sup$, the conditions
(BP1)---(BP3) of \cite{En} hold for the family
$\mathscr{B}^\star(x)$ and $\mathscr{B}^\star(x)$ is a base of a
Hausdorff topology $\tau^*$ at the point $x\in E^*$.

Let $c\in E$. Then there exists $b\in E$ such that $c<b<y$ and for
any open neighbourhood $V_b(y)$ of $y$ and any open neighbourhood
$V(c)$ such that $V(c)\subseteq E\setminus{\uparrow}b$ we have
$V_b(y)\cdot V(c)\subseteq V(c)$. We
also
have $V_b(y)\cdot
V_b(y)\subseteq V_b(y)$ for all $b\in{\downarrow}A=E$. Therefore
$(E^\star,\tau^\star)$ is a Hausdorff topological semilattice
which contains $E$ as a dense non-closed subsemilattice. This
contradicts the assumption that $E$ is an $H$-closed
topological semilattice.
The obtained contradictions imply that every subset of the
semilattice $E$ has  $\sup$.
The proof of the fact that every subset of $E$ has an $\inf$ is
similar.

Next we show that for every $H$-closed linearly ordered
topological semilattice $E$ condition $(ii)$ holds. Suppose that
there exists $x\in E$ such that
$x=\sup\left({\downarrow}x\setminus\{ x\}\right)$ and
$x\notin\operatorname{cl}_E\left({\downarrow}x\setminus\{
x\}\right)$. Since the topological semilattice $E$ is linearly
ordered, $L^\circ(x)={\downarrow}x\setminus\{ x\}$ is a clopen
subset of $E$.

Let $g\notin E$. We extend the semilattice operation from $E$ onto
$E^\circ=E\cup\{ g\}$ as follows:
\begin{equation*}
    g\cdot s=s\cdot g=
    \left\{%
\begin{array}{ll}
    g, & \hbox{if~} s\in{\uparrow}x; \\
    s, & \hbox{if~} s\in L^\circ(x).\\
\end{array}%
\right.
\end{equation*}
Obviously, the semilattice operation on $E^\circ$ determines a
linear order on $E^\circ$.

We define a topology $\tau^\circ$ on $E^\circ$ as follows. Let
$\tau$ be the topology on $E$. At any point $a\in
E=E^\circ\setminus\{ g\}$ bases of topologies $\tau^\circ$ and
$\tau$ coincide. We put
\begin{equation*}
 \mathscr{B}^\circ(g)=\left\{ U_s(g)=\{
 g\}\cup L^\circ(x)\setminus{\downarrow}s \mid
 s\in L^\circ(x)\right\}.
\end{equation*}
Since $\sup L^\circ(x)=x$, the set $U_s(g)$ is non-singleton for
any $s\in L^\circ(x)$. Therefore the  conditions (BP1)---(BP3) of
\cite{En} hold for the family $\mathscr{B}^\circ(g)$ and
$\mathscr{B}^\circ(g)$ is a base of the topology $\tau^\circ$ at
the point $g\in E^\circ$. Also since the set ${\uparrow}x$ is
closed in $(S^\circ,\tau^\circ)$ and $\sup L^\circ(x)=x$, we have
that the topology $\tau^\circ$ is Hausdorff. The proof of the
continuity of the semilattice operation in $(S^\circ,\tau^\circ)$
is similar as for
$(S^*,\tau^*)$ and $(S^\star,\tau^\star)$. Thus
condition $(ii)$ holds.

The proof of the assertion that if $E$ is a linearly ordered
$H$-closed topological semilattice, then the condition $(iii)$
holds, is similar.
Therefore the proof of the theorem is complete.
\end{proof}

Since the conditions $(i)$---$(iii)$ of Theorem~\ref{theorem1}
are preserved
by continuous homomorphisms, we have the following:

\begin{theorem}\label{theorem2}
Every linearly ordered $H$-closed topological semilattice is
absolutely $H$-closed.
\end{theorem}

Theorem~\ref{theorem1}
also
implies the following:

\begin{corollary}\label{corollary3}
Every linearly ordered $H$-closed topological semilattice contains
maximal and minimal idempotents.
\end{corollary}

Theorems~\ref{theorem1} and \ref{theorem2} imply
the following:

\begin{corollary}\label{corollary4}
Let $E$ be a linearly ordered $H$-closed topological semilattice
and $e\in E$. Then ${\uparrow}e$ and ${\downarrow}e$ are
(absolutely) $H$-closed topological semilattices.
\end{corollary}

\begin{theorem}\label{theorem6}
Every linearly ordered topological semilattice is a dense
subsemilattice of an $H$-closed linearly ordered topological
semilattice.
\end{theorem}

\begin{proof}
Let $E$ be a linearly ordered topological semilattice and let
$E_a$ be an algebraic copy of $E$. Then $E_a$ with operation
$\inf$ and $\sup$ is a lattice. It is well known that every
lattice embeds into a complete lattice (cf.
\cite[Theorem~V.2.1]{Birkhoff}). In our construction we shall use
the idea of proofs of Theorem~V.2.1 and Lemma~V.2.1 in
\cite{Birkhoff}.  We denote the lattice of all
ideals of $E_a$
by $\widetilde{E}_a$. Then pointwise operations $\inf$ and $\sup$ on
$\widetilde{E}_a$ coincide with $\bigcap$ and $\bigcup$ on
$\widetilde{E}_a$, respectively. Since $E_a$ is a linearly ordered
semilattice, we can identify $E_a$ with the subsemilattice of all
principal ideals of $E_a$ in $\widetilde{E}_a$.

For an ideal $I\in\widetilde{E}_a$ and a principal ideal
$I_e\in\widetilde{E}_a$ generated by an idempotent $e$ we put
$I_e\rho I$ if and only if every open neighbourhood of $e$
intersects $I$. Since $E_a$ and $\widetilde{E}_a$ are linearly
ordered semilattices, for principal ideals $I_e$ and $I_f$
generated by idempotents $e$ and $f$ from $E_a$, respectively, we
have $I_e\rho I_f$ if and only if $I_e=I_f$, i.e. $e=f$. We put
$\alpha=\Delta\cap\rho\cup\rho^{-1}$. Obviously the relation
$\alpha$ is an equivalence on $\widetilde{E}_a$. Since
$\widetilde{E}_a$ is a linearly ordered semilattice, $\alpha$ is a
congruence on $\widetilde{E}_a$ and hence
$\widetilde{E}=\widetilde{E}_a/\alpha$ is also
a linearly ordered
semilattice. We observe that $E_a$ is a subsemilattice of
$\widetilde{E}$.

We define a topology $\widetilde{\tau}$ on $\widetilde{E}$ as
follows. Let $\tau$ be the topology on $E$. At any point $a\in
E_a\subset\widetilde{E}$ bases of topologies $\widetilde{\tau}$
and $\tau$ coincide. For $x\in\widetilde{E}\setminus E_a$ we put
\begin{equation*}
 \widetilde{\mathscr{B}}(x)=\left\{ V_b(x)={\downarrow}x\setminus{\downarrow}b\mid
 b\in{\downarrow}x\cap E_a\right\}.
\end{equation*}
Then conditions (BP1)---(BP3) of \cite{En} hold for the family
$\widetilde{\mathscr{B}}(x)$ and $\widetilde{\mathscr{B}}(x)$ is a
base of a topology $\widetilde{\tau}$ at the point $x$ and since
$\tau$ is Hausdorff, so is $\widetilde{\tau}$. We
also observe
that the definition of $\widetilde{\tau}$ implies that
${\downarrow}e$ and ${\uparrow}e$ are closed subset of the
topological space $(\widetilde{E},\widetilde{\tau})$. Obviously
the semilattice operation on $(\widetilde{E},\widetilde{\tau})$ is
continuous, $(\widetilde{E},\widetilde{\tau})$ satisfies the
conditions $(i)$ and $(ii)$ of Theorem~\ref{theorem1} and $E$ is a
dense subsemilattice of $(\widetilde{E},\widetilde{\tau})$.

 We denote the lattice of all filters of
$\widetilde{E}$
by $\mathscr{F}(E)$. Then pointwise operations $\inf$ and $\sup$ on
$\mathscr{F}(E)$ coincide with $\bigcup$ and $\bigcap$ on
$\mathscr{F}(E)$, respectively. Since $\widetilde{E}$ is a
linearly ordered semilattice, we can identify $\widetilde{E}$ with
the subsemilattice in $\mathscr{F}(E)$ of all principal filters of
$\widetilde{E}$. By the dual theorem to Theorem~V.2.1 of
\cite{Birkhoff} the lattice $\mathscr{F}(E)$ is complete, and
since $\widetilde{E}$ is linearly ordered, so is $\mathscr{F}(E)$.

For a filter $F\in\widetilde{E}_a$ and a principal filter
$F_e\in\mathscr{F}(E)$ generated by an idempotent
$e\in\widetilde{E}$ we put $F_e\widetilde{\rho} F$ if and only if
every open neighbourhood of $e$ intersects $F$. Since
$\widetilde{E}$ and $\mathscr{F}(E)$ are linearly ordered
semilattices, for principal filters $F_e$ and $F_f$ generated by
idempotents $e$ and $f$ from $\widetilde{E}$, respectively, we
have $F_e\widetilde{\rho} F_f$ if and only if $F_e=F_f$, i.e.
$e=f$. We put
$\widetilde{\alpha}=\Delta\cap\widetilde{\rho}\cup{\widetilde{\rho}}\,^{-1}$.
Obviously $\widetilde{\alpha}$ is an equivalence on
$\mathscr{F}(E)$. Since $\mathscr{F}(E)$ is a linearly ordered
semilattice, $\widetilde{\alpha}$ is a congruence on
$\mathscr{F}(E)$ and hence
$\widetilde{\mathscr{F}}(E)=\mathscr{F}(E)/\widetilde{\alpha}$ is
a linearly ordered semilattice.

We define a topology $\tau_{\mathscr{F}}$ on
$\widetilde{\mathscr{F}}(E)$ as follows. At any point
$a\in\widetilde{E}\subseteq\widetilde{\mathscr{F}}(E)$ bases of
topologies $\tau_{\mathscr{F}}$ and $\widetilde{\tau}$ coincide.
For $x\in\widetilde{\mathscr{F}}(E)\setminus \widetilde{E}$ we put
\begin{equation*}
 \widetilde{\mathscr{B}}_{\mathscr{F}}(x)=\{
 W^b(x)={\uparrow}x\setminus{\uparrow}b\mid
 b\in{\uparrow}x\cap \widetilde{E}\}.
\end{equation*}
Then conditions (BP1)---(BP3) of \cite{En} hold for the family
$\widetilde{\mathscr{B}}_{\mathscr{F}}(x)$ and
$\widetilde{\mathscr{B}}_{\mathscr{F}}(x)$ is a base of a topology
$\tau_{\mathscr{F}}$ at the point $x$ and since $\widetilde{\tau}$
is Hausdorff, so is ${\tau}_{\mathscr{F}}$. Also we observe that
the definition of ${\tau}_{\mathscr{F}}$ implies that
${\downarrow}e$ and ${\uparrow}e$ are closed subsets of the
topological space
$(\widetilde{\mathscr{F}}(E),\tau_{\mathscr{F}})$, and hence the
semilattice operation on
$(\widetilde{\mathscr{F}}(E),\tau_{\mathscr{F}})$ is continuous,
$(\widetilde{E},\widetilde{\tau})$ satisfies the conditions $(i)$
and $(iii)$ of Theorem~\ref{theorem1} and $E$ is a dense
subsemilattice of
$(\widetilde{\mathscr{F}}(E),\tau_{\mathscr{F}})$.

Further we shall show that the condition $(ii)$ of
Theorem~\ref{theorem1} holds for the topological semilattice
$(\widetilde{\mathscr{F}}(E),\tau_{\mathscr{F}})$. Suppose to the
contrary that there exists a lower subset $A$ of $\mathscr{F}(E)$
such that $\sup A=x\notin A$ and
$x\notin\operatorname{cl}_{\mathscr{F}(E)} A$. Then
$A=\mathscr{F}(E)\setminus{\uparrow}x$ and $A$ are clopen subsets
of the topological space $(\mathscr{F}(E),\tau_{\mathscr{F}})$.
Since $x=\sup A\notin A$, there exists an increasing family of
ideals $\mathscr{I}=\{ I_\alpha\mid\alpha\in\mathcal{A}\}$ of the
semilattice $E_a$ such that $I_\alpha\subset I_\beta$ whenever
$\alpha<\beta$, $\alpha,\beta\in\mathcal{A}$,
$\sup\bigcup\mathscr{I}=x$, and $\bigcup\mathscr{I}\subset A$. The
existence of the family $\mathscr{I}$ follows from the fact that
the semilattice $\widetilde{E}$ is a dense subsemilattice of
$(\widetilde{\mathscr{F}}(E),\tau_{\mathscr{F}})$. However,
$\bigcup\mathscr{I}$ is an ideal in $E_a$ and hence
$\sup\bigcup\mathscr{I}\in A$, a contradiction. The obtained
contradiction implies  that statement $(ii)$ holds for the
topological semilattice
$(\widetilde{\mathscr{F}}(E),\tau_{\mathscr{F}})$.
\end{proof}

\begin{example}\label{example1}
Let ${\mathbb N}$ be the set of positive integers. Let $\{ x_n\}$
be an increasing sequence in ${\mathbb N}$. Put ${\mathbb
N}^\ast=\{ 0\}\cup\{\frac{1}{n}\mid n\in{\mathbb N}\}$. We define
the semilattice operation on ${\mathbb N}^\ast$ as follows
$ab=\min\{ a, b\}$, for $a, b\in{\mathbb N}^\ast$. Obviously, $0$
is the zero
element
of ${\mathbb N}^\ast$. We put $U_n(0)=\{
0\}\cup\{\frac{1}{x_k}\mid k\geqslant n\}$, $n\in{\mathbb N}$. A
topology $\tau$ on ${\mathbb N}^\ast$ is defined as follows: all
nonzero elements of ${\mathbb N}^\ast$ are isolated points in
${\mathbb N}^\ast$ and ${\mathscr B}(0)=\{ U_n(0)\mid n\in{\mathbb
N}\}$ is the base of the topology $\tau$ at the point
$0\in{\mathbb N}^\ast$. It is easy to see that $({\mathbb
N}^\ast,\tau)$ is a countable linearly ordered $\sigma$-compact
$0$-dimensional scattered locally compact metrizable topological
semilattice and if $x_{k+1}>x_k+1$ for every $k\in{\mathbb N}$,
then $({\mathbb N}^\ast,\tau)$ is a non-compact semilattice. We
also
observe that the family $\operatorname{Hom}(E,\{ 0,1\})$ of all
homomorphisms from a topological semilattice $E$ into the discrete
semilattice $(\{ 0,1\},\min)$ separates points for the topological
semilattice $E$.
\end{example}

Theorem~\ref{theorem1} implies the following:

\begin{proposition}\label{proposition1}
$({\mathbb N}^\ast,\tau)$ is an $H$-closed topological
semilattice.
\end{proposition}

\begin{remark}\label{remark10}
Example~\ref{example1} implies negative answers to the
following questions:
\begin{itemize}
\item[(1)] Is every closed subsemilattice of an $H$-closed topological
semilattice $H$-closed?
\item[(2)] (I.~Guran)
Does every locally compact topological semilattice embed into
a compact semilattice? 
\item[(3)] (cf. \cite{BH})
Does every globally bounded topological inverse Clifford semigroup embed into
a compact semigroup?
\item[(4)] Does every locally compact topological semilattice have a
base with open order convex subsets?
\end{itemize}
\end{remark}

\begin{remark}\label{remark11} Theorem~\ref{theorem2} and Example~\ref{example1}
imply that a closed subsemilattice of  an absolutely $H$-closed
topological semilattice is not $H$-closed.
\end{remark}

Example~\ref{example12} implies that there exist topologically
isomorphic linearly ordered topological semilattices $E_1$ and
$E_2$ which is dense subsemilattice of linearly ordered
topological semilattices $S_1$ and $S_2$, respectively, such that
$S_1$ and $S_2$ are not algebraically isomorphic.

\begin{example}\label{example12}
Let ${\mathbb N}$ be the set of positive integers. Let
$S_1=\{-\frac{1}{n}\mid n\in{\mathbb N}\}\cup\{
0\}\cup\{\frac{1}{n}\mid n\in{\mathbb N}\}$ and
$S_2=\{-1-\frac{1}{n}\mid n\in{\mathbb N}\}\cup\{-1\}\cup\{
0\}\cup\{\frac{1}{n}\mid n\in{\mathbb N}\}$ with usual topology
and operation $\min$. Then $E_1=\{-\frac{1}{n}\mid n\in{\mathbb
N}\}\cup\{\frac{1}{n}\mid n\in{\mathbb N}\}$ and
$E_2=\{-1-\frac{1}{n}\mid n\in{\mathbb N}\}\cup\{\frac{1}{n}\mid
n\in{\mathbb N}\}$ is discrete isomorphic semilattice, but the
semilattices $S_1$ and $S_2$ are not algebraically isomorphic.
\end{example}

Theorem~\ref{theorem13} gives a method of constructing new
$H$-closed and absolutely $H$-closed topological semilattices from
old.

\begin{theorem}\label{theorem13} Let
$S=\bigcup_{\alpha\in{\mathscr{A}}}S_{\alpha}$ be a topological
semilattice such that:
\begin{itemize}
    \item[$(i)$] $S_{\alpha}$ is an (absolutely) $H$-closed
    topological semilattice for any $\alpha\in{\mathscr{A}}$; and
    \item[$(ii)$] there exists an (absolutely) $H$-closed
    topological semilattice $T$ such that $T\subseteq S$ and
    $S_{\alpha}S_{\beta}\subseteq T$ for all $\alpha\neq\beta$,
    $\alpha,\beta\in{\mathscr{A}}$.
\end{itemize}
Then $S$ is an (absolutely) $H$-closed topological semilattice.
\end{theorem}

\begin{proof} We shall consider only the case when $S$ is an absolutely
$H$-closed topological semilattice. The proof in the other case is
similar.
Let $h\colon S\rightarrow G$ be a continuous homomorphism from $S$
into a topological semilattice $G$. Without loss of generality we
may
assume that $\operatorname{cl}_{G}(h(S))=G$.

Suppose that $G\setminus h(S)\neq\varnothing$. We fix $x\in
G\setminus h(S)$. The absolute $H$-closedness of the topological
semilattice $T$ implies that there exists an open neighbourhood
$U(x)$ of the point $x$ in $G$ such that $U(x)\cap
h(T)=\varnothing$. Since $G$ is a topological semilattice, there
exists an open neighbourhood $V(x)$ of $x$ in $G$ such that
$V(x)V(x)\subseteq U(x)$. Since the topological semilattice
$S_{\alpha}$ is absolutely $H$-closed, the neighbourhood $V(x)$
intersects infinitely many semilattices $h(S_{\beta})$,
$\beta\in{\mathscr{A}}$. Therefore $V(x)V(x)\cap
h(T)\neq\varnothing$. This is in disagreement
with
the choice of the
neighbourhood $U(x)$. This contradiction implies the
assertion
of the theorem.
\end{proof}


\section*{Acknowledgements}
This research was supported by SRA grants P1-0292-0101-04 and  BI-UA/07-08-001.
The authors thank the referee for important remarks
and  suggestions. 

\vfill\eject

\end{document}